\input amstex
\documentstyle{amsppt}

\document
\magnification=1200
\NoBlackBoxes
\nologo
\vsize18cm
\centerline{\bf Functional realization of some elliptic}
\medskip
\centerline{\bf Hamiltonian structures and bosonization}
\medskip
\centerline{\bf of the corresponding quantum algebras}
\bigskip
\centerline{\bf by A.V.Odesskii and B.L.Feigin}
\bigskip
 \centerline{\bf Introduction}
\medskip

Let $\frak{p}$ be the parabolic subalgebra of some semisimple Lie algebra $\frak{g}$ and $P$ the corresponding group. Let $\Cal M(\Cal E,\frak{p})$ be the moduli space of $P$-bundles on the elliptic curve $\Cal E$. In [5] we define the Hamiltonian structure on the manifold $\Cal M(\Cal E,\frak{p})$. There is a natural problem: to quantize the coordinate ring of each connected component of $\Cal M(\Cal E,\frak{p})$.

We denote by $Q_n(\Cal E,\tau)$ the corresponding quantum algebras in the case $\frak{g}=sl_2$. Here $\tau\in\Cal E$ is a parameter of quantization, $n\in\Bbb N$ is a number of the connected component of $\Cal M(\Cal E,\frak{p})$. This component is isomorphic to $\Bbb P^{n-1}$ in this case. So its coordinate ring is isomorphic to the polynomial ring in $n$ variables and the algebra $Q_n(\Cal E,\tau)$ is a graded deformation of this polynomial ring. We denote the corresponding Poisson algebra by $q_n(\Cal E)$. 

More generally, we denote by $Q_{n,k}(\Cal E,\tau)$ the corresponding quantum algebras in the case $\frak{g}=sl_{k+1}$ and $\frak{p}$ is a parabolic subalgebra for the flag $0\subset V\subset\Bbb C^{k+1}$, $\roman{dim}V=1$.  Here $\tau\in\Cal E$ is a parameter of quantization, $n\in\Bbb N$ is a number of the connected component of $\Cal M(\Cal E,\frak{p})$. If $n$ and $k$ have no common divisors, then this component is isomorphic to $\Bbb P^{n-1}$. So its coordinate ring is isomorphic to the polynomial ring in $n$ variables and the algebra $Q_{n,k}(\Cal E,\tau)$ is a graded deformation of this polynomial ring. We denote the corresponding Poisson algebra by $q_{n,k}(\Cal E)$.  We have $Q_n(\Cal E,\tau)=Q_{n,1}(\Cal E,\tau)$, $q_n(\Cal E)=q_{n,1}(\Cal E)$.

In the papers [1,2] we constructed the family of associative algebras $Q_n(\Cal E,\tau)$. The algebra $Q_n(\Cal E,\tau)$ is $\Bbb Z$-graded and depends on 2 continuous parameters: an elliptic curve $\Cal E=\Bbb C\diagup\Gamma$ and a point $\tau\in\Cal E$. We have $Q_n(\Cal E,\tau)=\Bbb C\oplus F_1\oplus F_2\oplus\dots$ and $F_\alpha*F_\beta\subset F_{\alpha+\beta}$. The Hilbert function is $1+\botshave{\sum_{\alpha\geqslant 1}}\roman{dim}\,F_\alpha t^\alpha=(1-t)^{-n}$. If $\tau=0$ then the product $*$ is commutative and the algebra $Q_n(\Cal E,0)$ is a polynomial ring in $n$ variables. So the algebra $Q_n(\Cal E,0)$ does not depend on $\Cal E$. Considering the product $*$ in a neighborhood of $\tau=0$ with the fixed curve $\Cal E$ we will have a Poisson structure on the polynomial ring $Q_n(\Cal E,0)$ that depends on $\Cal E$. We call the following construction of the algebra $Q_n(\Cal E,\tau)$ its functional realization, because the graded components $F_\alpha$ are described as the spaces of functions and the product $*$ is given by an explicit formula. Namely, $F_\alpha$ is a space of holomorphic symmetric functions in $\alpha$ variables $f(z_1,\dots,z_\alpha)$ with the properties: $f(z_1+1,z_2,\dots,z_\alpha)=f(z_1,\dots,z_\alpha), f(z_1+\eta,z_2,\dots,z_\alpha)=e^{-2\pi i(nz_1-\alpha n\tau)}f(z_1,\dots,z_\alpha)$. Here $1,\eta$ are generators of the lattice $\Gamma$, $\roman{Im}\,\eta>0$. It is clear that $F_\alpha\cong S^\alpha(\Theta_{n,-\alpha n\tau}(\Gamma))$ (see {\bf Notations}) and $\roman{dim}\,F_\alpha=\frac{n(n+1)\dots(n+\alpha-1)}{\alpha!}$. The product $*$ is given in the following way: for $f\in F_\alpha$, $g\in F_\beta$ we have
 
$$f*g(z_1,\dots,z_{\alpha+\beta})=\eqno(1)$$

$$\frac{1}{\alpha!\beta!}\botshave{\sum_{\sigma\in S_{\alpha+\beta}}}f(z_{\sigma_1},\dots,z_{\sigma_\alpha})g(z_{\sigma_{\alpha+1}}-2\alpha\tau,\dots,z_{\sigma_{\alpha+\beta}}-2\alpha\tau)\prod\Sb1\leqslant p\leqslant\alpha\\ \\\alpha+1\leqslant q\leqslant\alpha+\beta\endSb \frac{\theta(z_{\sigma_p}-z_{\sigma_q}-n\tau)}{\theta(z_{\sigma_p}-z_{\sigma_q})}$$
If $\tau=0$ then the formula (1) gives the usual product in the symmetric algebra $S^*(\Theta_{n,0}(\Gamma))$. Decompose the right side of the formula (1) in the Taylor series $f*g=fg+c_1(f,g)\tau+o(\tau)$, where $fg$ is the product in the symmetric algebra $S^*(\Theta_{n,0}(\Gamma))$. Then we will have the explicit formula for the Poisson structure on the algebra $S^*(\Theta_{n,0}(\Gamma))$. As usual, $\{f,g\}=c_1(f,g)-c_1(g,f)$. In the explicit form:

$$\{f,g\}(z_1,\dots,z_{\alpha+\beta})=$$

$$\frac{1}{\alpha!\beta!}\botshave{\sum_{\sigma\in S_{\alpha+\beta}}}\bigg(-2nf(z_{\sigma_1},\dots,z_{\sigma_\alpha})g(z_{\sigma_{\alpha+1}},\dots,z_{\sigma_{\alpha+\beta}})\sum\Sb1\leqslant p\leqslant\alpha\\ \\\alpha+1\leqslant q\leqslant\alpha+\beta\endSb\frac{\theta^\prime(z_{\sigma_p}-z_{\sigma_q})}{\theta(z_{\sigma_p}-z_{\sigma_q})}+$$

$$+2\beta g(z_{\sigma_{\alpha+1}},\dots,z_{\sigma_{\alpha+\beta}})\botshave{\sum_{1\leqslant p\leqslant\alpha}}f^\prime_{z_{\sigma_p}}(z_{\sigma_1},\dots,z_{\sigma_\alpha})-$$

$$-2\alpha f(z_{\sigma_1},\dots,z_{\sigma_\alpha})\botshave{\sum_{\alpha+1\leqslant q\leqslant\alpha+\beta}}g^\prime_{z_{\sigma_q}}(z_{\sigma_{\alpha+1}},\dots,z_{\sigma_{\alpha+\beta}})\bigg)$$

There is another construction of the algebra $Q_n(\Cal E,\tau)$ that we call bosonization. Let $A_p(\Cal E,\tau)$ be the algebra generated by $\{e_1,\dots,e_p,\varphi(z_1,\dots,z_p)\}$, where $\varphi$ is any meromorphic functions in variables $z_1,\dots,z_p$. We assume that the following relations hold: 

$$z_\alpha z_\beta=z_\beta z_\alpha, e_\alpha z_\beta=(z_\beta-2\tau)e_\alpha, e_\alpha z_\alpha=(z_\alpha+(n-2)\tau)e_\alpha$$

$$e_\alpha e_\beta=-e^{2\pi i(z_\beta-z_\alpha)}\frac{\theta(z_\alpha-z_\beta-n\tau)}{\theta(z_\beta-z_\alpha-n\tau)}e_\beta e_\alpha,\text{ here } \alpha\ne\beta$$

There is a homomorphism $x: Q_n(\Cal E,\tau)\to A_p(\Cal E,\tau)$ such that for $f\in F_1=\Theta_{n,-n\tau}(\Gamma)$ we have: $x(f)=\botshave{\sum_{1\leqslant\alpha\leqslant p}}f(z_\alpha)e_\alpha$.
It is clear that for $\tau=0$ the algebra $A_p(\Cal E,0)$ is commutative and $x$ is the homomorphism of the Poisson algebras. It is easy to see that the functional realization of the algebra $Q_n(\Cal E,\tau)$ (and the corresponding Poisson algebra) follows from the bosonization.

In [1,2] we studied the algebras $Q_{n,k}(\Cal E,\tau)$. Here $1\leqslant k<n$; $n$ and $k$ have no common divisors. We have $Q_n(\Cal E,\tau)=Q_{n,1}(\Cal E,\tau)$. In the general case the algebra $Q_{n,k}(\Cal E,\tau)$ is also $\Bbb Z$-graded, so we have: $Q_{n,k}(\Cal E,\tau)=\Bbb C\oplus P_1\oplus P_2\oplus\dots$ and $P_\alpha*P_\beta\subset P_{\alpha+\beta}$. The Hilbert function is $1+\botshave{\sum_{\alpha\geqslant1}}\roman{dim}\,P_\alpha t^\alpha=(1-t)^{-n}$. If $\tau=0$ than the algebra $Q_{n,k}(\Cal E,0)$ is the polynomial ring in $n$ variables. We denote the corresponding Poisson algebra by $q_{n,k}(\Cal E)$. Let $\frac{n}{k}=n_1-\frac{1}{n_2-\dots-\frac{1}{n_p}}$ be the decomposition into the continuous fraction, where $n_1,\dots,n_p\geqslant2$. It is clear that $n=d(n_1,\dots,n_p)$, $k=d(n_2,\dots,n_p)$. The space $P_1$ is isomorphic to the space $\Theta_{(n_1,\dots,n_p)}(\Gamma)$ (see {\bf Notations}). We have $\roman{dim}\,\Theta_{(n_1,\dots,n_p)}(\Gamma)=n$.

In [3,4] we consider the case when $\frak{g}$ is general and  $\frak{p}$ is a Borel subalgebra. We denote the corresponding quantum algebras by $Q_{n,\Delta}(\Cal E,\tau)$. Here $\Delta$ is the root system of $\frak{g}$, $n: L\to\Bbb Z$ is the homomorphism of the additive groups, $L$ is the lattice generated by $\Delta$.  

In this paper we consider the case $\frak{g}=sl_N$ and a general parabolic subalgebra $\frak{p}$. We introduce a functional realization of the Poisson algebras in this case. We also introduce a certain construction of the corresponding quantum algebras. We denote these quantum algebras by $Q_{n_1,k_1}\widehat\otimes\dots\widehat\otimes Q_{n_h,k_h}(\Cal E,\tau)$ and the corresponding Poisson algebras by $q_{n_1,k_1}\widehat\otimes\dots\widehat\otimes q_{n_h,k_h}(\Cal E)$. 

\medskip
Now we describe the contents of the paper.
\medskip

 In {\bf \S1} we construct the functional realization of the Poisson algebra $q_{n,k}(\Cal E)$. In this construction $P_\alpha=S^\alpha(\Theta_{(n_1,\dots,n_p)}(\Gamma))$ is realized as the space of holomorphic functions $f(x_{1,1},\dots,x_{p,1};\dots;x_{1,\alpha},\dots,x_{p,\alpha})$ satisfying some properties. The Poisson bracket is given by the formulas (2) and (3).
 
In {\bf \S2} we construct an analogue of bosonization of the algebras $Q_{n,k}(\Cal E,\tau)$.

In {\bf \S\S3-5} we construct the algebras that we denote by $Q_{n_1,k_1}\widehat\otimes\dots\widehat\otimes Q_{n_h,k_h}(\Cal E,\tau)$ and the corresponding Poisson algebras $q_{n_1,k_1}\widehat\otimes\dots\widehat\otimes q_{n_h,k_h}(\Cal E)$. These algebras are $\Bbb Z^h$-graded, so $Q_{n_1,k_1}\widehat\otimes\dots\widehat\otimes Q_{n_h,k_h}(\Cal E,\tau)=\botshave{\bigoplus_{\alpha_1,\dots,\alpha_h\geqslant0}}P_{\alpha_1,\dots,\alpha_h}$ and $P_{\alpha_1,\dots,\alpha_h}*P_{\beta_1,\dots,\beta_h}\subset P_{\alpha_1+\beta_1,\dots,\alpha_h+\beta_h}$. The space $M=P_{0,\dots,0}$ is the field of meromorphic functions in variables $z_{1,2},\dots,z_{h-1,h}$. The dimensions $\roman{dim}_MP_{\alpha_1,\dots,\alpha_h}$ are finite and the Hilbert function is $$\botshave{\sum_{\alpha_1,\dots,\alpha_h\geqslant0}}\roman{dim}_MP_{\alpha_1,\dots,\alpha_h}t_1^{\alpha_1}\dots t_h^{\alpha_h}=\botshave{\prod_{1\leqslant\lambda\leqslant\nu\leqslant h}}(1-t_\lambda t_{\lambda+1}\dots t_\nu)^{-d(N_\lambda\Delta N_{\lambda+1}\Delta\dots\Delta N_\nu)}$$
Here $N_\alpha=(n_{1,\alpha},n_{2,\alpha},\dots,n_{p_\alpha,\alpha})$, where $\frac{n_\alpha}{k_\alpha}=n_{1,\alpha}-\frac{1}{n_{2,\alpha}-\dots-\frac{1}{n_{p_\alpha,\alpha}}}$ is the decomposition into the continuous fraction, $1\leqslant\alpha\leqslant h$, $n_{i,j}\geqslant2$. See {\bf Notations} for the definition $\Delta$.

In {\bf \S3} we construct the functional realization of the Poisson algebra $q_{n,k}\widehat\otimes q_{m,l}(\Cal E)$.
In {\bf \S4} we construct the functional realization of the Poisson algebra $q_{n_1,k_1}\widehat\otimes\dots\widehat\otimes q_{n_h,k_h}(\Cal E)$ in general case.

In {\bf \S5} we construct the bosonization of the algebra $Q_{n_1,k_1}\widehat\otimes\dots\widehat\otimes Q_{n_h,k_h}(\Cal E,\tau)$ in general case.
\medskip
In [3,4] we consider the case $k_1=\dots=k_h=1$. In this case $\frak{g}=sl_{h+1}$ and $\frak{p}$ is a Borel subalgebra. In those papers the algebra $Q_{n_1,1}\widehat\otimes\dots\widehat\otimes Q_{n_h,1}(\Cal E,\tau)$ was denoted by $Q_{n,\Delta}(\Cal E,\tau)$. Here $\Delta$ is the root system $A_h$, $n: L\to\Bbb Z$ is the homomorphism of the additive groups, $L$ is the lattice generated by $\Delta$, $n_\alpha=n(\delta_\alpha)$ where $\delta_1,\dots,\delta_h$ are simple positive roots. In [3] \S1 and [4] \S1 we constructed the functional realization of this algebra. In [3] \S1 and [4] \S2 we constructed the bosonization of this algebra.

\newpage
\centerline{\bf Notations}
\medskip

Let $\Cal E=\Bbb C\diagup\Gamma$ be an elliptic curve, where $\Gamma=\{m_1+m_2\eta; m_1,m_2\in \Bbb Z\}$ is a lattice, $\roman{Im}\,\eta>0$. For $m\in \Bbb Z$, $c\in\Bbb C$ we denote by $\Theta_{m,c}(\Gamma)$ the space of holomorphic functions $f(z)$ with the following properties: $f(z+1)=f(z), f(z+\eta)=e^{-2\pi i(mz+c)}f(z)$. It is clear that $\roman{dim}\,\Theta_{m,c}(\Gamma)=m$ if $m>0$. For $m>0$ the elements of $\Theta_{m,c}(\Gamma)$ are called $\theta$-functions of order $m$. It is easy to check that every $\theta$-function of order $m$ has exactly $m$ zeros $\roman{mod}\,\Gamma$ and the sum of these zeros is equal to $c+\frac{1}{2}m$ $\roman{mod}\,\Gamma$. Let $\theta(z)=\botshave{\sum_{\alpha\in\Bbb Z}}(-1)^\alpha e^{2\pi i(\alpha z+\frac{\alpha(\alpha-1)}{2}\eta)}$. It is clear that $\theta(z)\in\Theta_{1,\frac{1}{2}}(\Gamma)$, $\theta(0)=0$, $\theta(-z)=-e^{-2\pi iz}\theta(z)$.
 
For the sequence of natural numbers $(n_1,\dots,n_p)$ we denote   $d(n_1,\dots,n_p)=\roman{det}\,(m_{i,j})$, here $(m_{i,j})$ is the $p\times p$ matrix with the elements $m_{i,i}=n_i, m_{i,i+1}=m_{i+1,i}=-1$, $m_{i,j}=0$ for $|i-j|>1$. For $p=0$ we assume $d(\emptyset)=1$. Let $n_i\geqslant 2$ for all $1\leqslant i\leqslant p$, $n=d(n_1,\dots,n_p)$, $k=d(n_2,\dots,n_p)$. It is clear that $n$ and $k$ have no common divisors, $1\leqslant k<n$ and $\frac{n}{k}=n_1-\frac{1}{n_2-\frac{1}{n_3-\dots-\frac{1}{n_p}}}$.

For two sequences $A=(a_1,\dots,a_p)$ and $B=(b_1,\dots,b_q)$ we denote $A\Delta B=(a_1,\dots,a_{p-1},a_p+b_1,b_2,\dots,b_q)$. It is a sequence of length $p+q-1$. Particularly, for $p=q=1$ we have $A\Delta B=(a_1+b_1)$. It is clear that the operation $\Delta$ is associative.

We denote by $\Theta_{(n_1,\dots,n_p)}(\Gamma)$ the space of holomorphic functions $f(z_1,\dots,z_p)$ with the following properties: 
$$f(z_1,\dots,z_\alpha+1,\dots,z_p)=f(z_1,\dots,z_p)$$
$$f(z_1,\dots,z_\alpha+\eta,\dots,z_p)=e^{-2\pi i(n_\alpha z_\alpha-z_{\alpha-1}-z_{\alpha+1})}f(z_1,\dots,z_p)$$
We assume here $z_0=z_{p+1}=0$. It is easy to check that $\roman{dim}\,\Theta_{(n_1,\dots,n_p)}(\Gamma)=d(n_1,\dots,n_p)$ for $n_1,\dots,n_p\geqslant 2$.

Let $N=(n_1,\dots,n_p)$, $M=(m_1,\dots,m_q)$ be two sequences of natural numbers. Let $f(z_1,\dots,z_p)\in \Theta_{(n_1,\dots,n_p)}(\Gamma)$, $g(z_1,\dots,z_q)\in  \Theta_{(m_1,\dots,m_q)}(\Gamma)$. It is clear that if $\varphi(z_1,\dots,z_{p+q-1})=f(z_1,\dots,z_p)g(z_p,\dots,z_{p+q-1})$ then $\varphi\in \Theta_{N\Delta M}(\Gamma)$.

Let $A$ be some commutative associative algebra. We remind that the structure of a Lie algebra $\{,\}$ on the space $A$ with the following property (Leibniz rule):
$\{f,gh\}=\{f,g\}h+g\{f,h\}$ for $f,g,h\in A$ is called a Poisson structure on $A$. The algebra $A$ with a Poisson structure is called a Poisson algebra. Let $*_t$ be such a family of associative products on the space $A$ that holomorphically depends on $t\in U\subset\Bbb C$, here $U$ is an open subset, $0\in U$, and $f*_tg=fg+\frac{1}{2}\{f,g\}t+o(t)$ as $t\to0$. Let $A_t$ be the associative algebra with the product $*_t$. The family of associative algebras $A_t$ is called a quantization of the Poisson algebra $A$.

\newpage

\centerline{\bf \S1. Functional realization of the Poisson algebra $q_{n,k}(\Cal E)$}
\medskip
Let $q_{n,k}(\Cal E)=S^*(\Theta_{(n_1,\dots,n_p)}(\Gamma))$, where $\frac{n}{k}=n_1-\frac{1}{n_2-\dots-\frac{1}{n_p}}$ is the decomposition into the continuous fraction, $n_1,\dots,n_p\geqslant2$. So $q_{n,k}(\Cal E)$ is the polynomial ring in $n$ variables and the space $S^\alpha(\Theta_{(n_1,\dots,n_p)}(\Gamma))$ of elements of degree $\alpha$ is realized as the space of holomorphic functions $f(x_{1,1},\dots,x_{p,1};\dots;x_{1,\alpha},\dots,x_{p,\alpha})$ with the following properties: 

1. Symmetry. For each $\sigma\in S_\alpha$ we have $f(x_{1,\sigma_1},\dots,x_{p,\sigma_1};\dots;x_{1,\sigma_\alpha},\dots,x_{p,\sigma_\alpha})=f(x_{1,1},\dots,x_{p,1};\dots;x_{1,\alpha},\dots,x_{p,\alpha})$.

2. Periodicity and quasiperiodicity.
 $$f(x_{1,1},\dots,x_{\mu,\nu}+1,\dots,x_{p,\alpha})=f(x_{1,1},\dots,x_{p,\alpha})$$

$$f(x_{1,1},\dots,x_{\mu,\nu}+\eta,\dots,x_{p,\alpha})=e^{-2\pi i(n_\mu x_{\mu,\nu}-x_{\mu-1,\nu}-x_{\mu+1,\nu})}f(x_{1,1},\dots,x_{p,\alpha})$$
We assume here $x_{0,\nu}=x_{p+1,\nu}=0$.

The product on the space $q_{n,k}(\Cal E)$ is given by the usual formula:

 for $f\in S^\alpha(\Theta_{(n_1,\dots,n_p)}(\Gamma))$, $g\in S^\beta(\Theta_{(n_1,\dots,n_p)}(\Gamma))$ we have:

$$fg(x_{1,1},\dots,x_{p,1};\dots;x_{1,\alpha+\beta},\dots,x_{p,\alpha+\beta})=\frac{1}{\alpha!\beta!}\botshave{\sum_{\sigma\in S_{\alpha+\beta}}}$$

$$f(x_{1,\sigma_1},\dots,x_{p,\sigma_1};\dots;x_{1,\sigma_\alpha},\dots,x_{p,\sigma_\alpha})g(x_{1,\sigma_{\alpha+1}},\dots,x_{p,\sigma_{\alpha+1}};\dots;x_{1,\sigma_{\alpha+\beta}},\dots,x_{p,\sigma_{\alpha+\beta}})$$

We define the operation $\{,\}$ on the space $q_{n,k}(\Cal E)$ in the following way: if $f, g\in \Theta_{(n_1,\dots,n_p)}(\Gamma)$ then for $\{f,g\}\in S^2(\Theta_{(n_1,\dots,n_p)}(\Gamma))$ we have: 

$$\{f,g\}(x_1,\dots,x_p;y_1,\dots,y_p)= \eqno(2)$$

$$\botshave{\sum_{1\leqslant\alpha\leqslant p}}\frac{d(n_1,\dots,n_{\alpha-1})+d(n_{\alpha+1},\dots,n_p)}{d(n_1,\dots,n_p)}\bigg(g(x_1,\dots,x_p)f^\prime_{y_\alpha}(y_1,\dots,y_p)+$$

$$+g(y_1,\dots,y_p)f^\prime_{x_\alpha}(x_1,\dots,x_p)-f(x_1,\dots,x_p)g^\prime_{y_\alpha}(y_1,\dots,y_p)-f(y_1,\dots,y_p)g^\prime_{x_\alpha}(x_1,\dots,x_p)\bigg)+$$

$$+\bigg(\frac{\theta^\prime(y_1-x_1)}{\theta(y_1-x_1)}+\frac{\theta^\prime(y_p-x_p)}{\theta(y_p-x_p)}-2\pi i\bigg)\bigg(f(x_1,\dots,x_p)g(y_1,\dots,y_p)-g(x_1,\dots,x_p)f(y_1,\dots,y_p)\bigg)+$$

$$\theta^\prime(0)\botshave{\sum_{1\leqslant\alpha<p}}\frac{\theta(x_\alpha+y_{\alpha+1}-y_\alpha-x_{\alpha+1})}{\theta(x_\alpha-y_\alpha)\theta(y_{\alpha+1}-x_{\alpha+1})}\bigg(f(y_1,\dots,y_\alpha,x_{\alpha+1},\dots,x_p)g(x_1,\dots,x_\alpha,y_{\alpha+1},\dots,y_p)-$$

$$-g(y_1,\dots,y_\alpha,x_{\alpha+1},\dots,x_p)f(x_1,\dots,x_\alpha,y_{\alpha+1},\dots,y_p)\bigg)$$

In the general case, if $f\in S^\alpha(\Theta_{(n_1,\dots,n_p)}(\Gamma))$, $g\in S^\beta(\Theta_{(n_1,\dots,n_p)}(\Gamma))$ then for $\{f,g\}\in S^{\alpha+\beta}(\Theta_{(n_1,\dots,n_p)}(\Gamma))$ we have by definition: 

$$\{f,g\}(x_{1,1},\dots,x_{p,1};\dots;x_{1,\alpha+\beta},\dots,x_{p,\alpha+\beta})=\frac{1}{\alpha!\beta!}\botshave{\sum_{\sigma\in S_{\alpha+\beta}}}\eqno(3)$$

$$\bigg(\beta g(x_{1,\sigma_{\alpha+1}},\dots,x_{p,\sigma_{\alpha+\beta}})\sum\Sb1\leqslant\psi\leqslant p\\ \\1\leqslant\mu\leqslant\alpha\endSb\frac{d(n_1,\dots,n_{\psi-1})+d(n_{\psi+1},\dots,n_p)}{d(n_1,\dots,n_p)}f^\prime_{x_{\psi,\sigma_\mu}}(x_{1,\sigma_1},\dots,x_{p,\sigma_\alpha})-$$

$$\alpha f(x_{1,\sigma_1},\dots,x_{p,\sigma_\alpha})\sum\Sb1\leqslant\psi\leqslant p\\ \\\alpha+1\leqslant\mu\leqslant\alpha+\beta\endSb\frac{d(n_1,\dots,n_{\psi-1})+d(n_{\psi+1},\dots,n_p)}{d(n_1,\dots,n_p)}g^\prime_{x_{\psi,\sigma_\mu}}(x_{1,\sigma_{\alpha+1}},\dots,x_{p,\sigma_{\alpha+\beta}})+$$

$$+\bigg(\sum\Sb1\leqslant\mu\leqslant\alpha\\ \\\alpha+1\leqslant\mu^\prime\leqslant\alpha+\beta\endSb\frac{\theta^\prime(x_{1,\sigma_{\mu^\prime}}-x_{1,\sigma_\mu})}{\theta(x_{1,\sigma_{\mu^\prime}}-x_{1,\sigma_\mu})}+\frac{\theta^\prime(x_{p,\sigma_{\mu^\prime}}-x_{p,\sigma_\mu})}{\theta(x_{p,\sigma_{\mu^\prime}}-x_{p,\sigma_\mu})}-2\pi i\alpha\beta\bigg)\times$$

$$\times f(x_{1,\sigma_1},\dots,x_{p,\sigma_\alpha})g(x_{1,\sigma_{\alpha+1}},\dots,x_{p,\sigma_{\alpha+\beta}})+$$

$$+\theta^\prime(0)\sum\Sb1\leqslant\mu\leqslant\alpha\\ \\\alpha+1\leqslant\mu^\prime\leqslant\alpha+\beta\\ \\1\leqslant\psi\leqslant p-1\endSb\frac{\theta(x_{\psi,\sigma_\mu}+x_{\psi+1,\sigma_{\mu^\prime}}-x_{\psi,\sigma_{\mu^\prime}}-x_{\psi+1,\sigma_\mu})}{\theta(x_{\psi,\sigma_\mu}-x_{\psi,\sigma_{\mu^\prime}})\theta(x_{\psi+1,\sigma_{\mu^\prime}}-x_{\psi+1,\sigma_\mu})}\times$$

$$\times f(x_{1,\sigma_1},\dots,x_{p,\sigma_1};\dots;x_{1,\sigma_{\mu^\prime}},\dots,x_{\psi,\sigma_{\mu^\prime}},x_{\psi+1,\sigma_\mu},\dots,x_{p,\sigma_\mu};\dots;x_{1,\sigma_\alpha},\dots,x_{p,\sigma_\alpha})\times$$

$$\times g(x_{1,\sigma_{\alpha+1}},\dots,x_{p,\sigma_{\alpha+1}};\dots;x_{1,\sigma_\mu},\dots,x_{\psi,\sigma_\mu},x_{\psi+1,\sigma_{\mu^\prime}},\dots,x_{p,\sigma_{\mu^\prime}};\dots;x_{1,\sigma_{\alpha+\beta}},\dots,x_{p,\sigma_{\alpha+\beta}})\bigg)$$

{\bf Proposition 1.} {\it The operation $\{,\}$ defines a Poisson bracket on the space $q_{n,k}(\Cal E)$.}

\newpage

\centerline{\bf \S2. Bosonization of the algebra $Q_{n,k}(\Cal E,\tau)$}
\medskip

Let $\frac{n}{k}=n_1-\frac{1}{n_2-\dots-\frac{1}{n_p}}$,  $n_1,\dots,n_p\geqslant2$. Let $A_{m_1,\dots,m_p}(\Cal E,\tau)$ be the associative algebra generated by $\{e_{\alpha_1,\dots,\alpha_p}; 1\leqslant\alpha_t\leqslant m_t, 1\leqslant t\leqslant p\}$ and $\{\varphi(y_{1,1},\dots,y_{p,m_p})\}$, where $\varphi$ is any meromorphic function in variables $\{y_{\lambda,\alpha}; 1\leqslant\lambda\leqslant p, 1\leqslant\alpha\leqslant m_\lambda\}$. We assume that the following relations hold (see also [2] for the cases $p=1$ and $p=2$):  

$$y_{\lambda,\alpha}y_{\nu,\beta}=y_{\nu,\beta}y_{\lambda,\alpha},$$

$$e_{\alpha_1,\dots,\alpha_p}y_{\nu,\beta}=(y_{\nu,\beta}-(d(n_1,\dots,n_{\nu-1})+d(n_{\nu+1},\dots,n_p))\tau)e_{\alpha_1,\dots,\alpha_p}, \text{ here } \alpha_\nu\ne\beta$$

$$e_{\alpha_1,\dots,\alpha_p}y_{\nu,\alpha_\nu}=(y_{\nu,\alpha_\nu}+(d(n_1,\dots,n_p)-d(n_1,\dots,n_{\nu-1})-d(n_{\nu+1},\dots,n_p))\tau)e_{\alpha_1,\dots,\alpha_p}$$

$$e_{\alpha_1,\dots,\alpha_p}e_{\beta_1,\dots,\beta_p}=\Lambda e_{\beta_1,\dots,\beta_p}e_{\alpha_1,\dots,\alpha_p}+\eqno(4)$$

$$+\botshave{\sum_{1\leqslant t\leqslant p-1}}\Lambda_{t,t+1}e_{\beta_1,\dots,\beta_t,\alpha_{t+1},\dots,\alpha_p}e_{\alpha_1,\dots,\alpha_t,\beta_{t+1},\dots,\beta_p},$$

Here $\alpha_1\ne\beta_1,\dots,\alpha_p\ne\beta_p$ and

$$\Lambda=\frac{e^{-2\pi in\tau}\theta(y_{1,\beta_1}-y_{1,\alpha_1})\theta(y_{p,\beta_p}-y_{p,\alpha_p}+n\tau)}{\theta(y_{1,\beta_1}-y_{1,\alpha_1}-n\tau)\theta(y_{p,\beta_p}-y_{p,\alpha_p})}\eqno(5)$$

$$\Lambda_{t,t+1}=\frac{e^{-2\pi in\tau}\theta(n\tau)\theta(y_{1,\beta_1}-y_{1,\alpha_1})}{\theta(y_{1,\beta_1}-y_{1,\alpha_1}-n\tau)}\cdot\frac{\theta(y_{t,\beta_t}+y_{t+1,\beta_{t+1}}-y_{t,\alpha_t}-y_{t+1,\alpha_{t+1}})}{\theta(y_{t,\beta_t}-y_{t,\alpha_t})\theta(y_{t+1,\beta_{t+1}}-y_{t+1,\alpha_{t+1}})}$$

In the general case, if some indexes are the same, then we have:

$$e_{\mu_1,\dots,\mu_{\psi-1},\mu_\psi,\alpha_1,\dots,\alpha_\varphi,\gamma_1,\gamma_2,\dots,\gamma_q}e_{\mu_1^\prime,\dots,\mu_{\psi-1}^\prime,\mu_\psi,\beta_1,\dots,\beta_\varphi,\gamma_1,\gamma_2^\prime,\dots,\gamma_q^\prime}=$$

$$\Lambda e_{\mu_1,\dots,\mu_{\psi-1},\mu_\psi,\beta_1,\dots,\beta_\varphi,\gamma_1,\gamma_2,\dots,\gamma_q}e_{\mu_1^\prime,\dots,\mu_{\psi-1}^\prime,\mu_\psi,\alpha_1,\dots,\alpha_\varphi,\gamma_1,\gamma_2^\prime,\dots,\gamma_q^\prime}+\eqno(6)$$

$$\botshave{\sum_{1\leqslant t<\varphi}}\Lambda_{t,t+1}e_{\mu_1,\dots,\mu_{\psi-1},\mu_\psi,\beta_1,\dots,\beta_t,\alpha_{t+1},\dots,\alpha_\varphi\gamma_1,\gamma_2,\dots,\gamma_q}e_{\mu_1^\prime,\dots,\mu_{\psi-1}^\prime,\mu_\psi,\alpha_1,\dots,\alpha_t,\beta_{t+1},\dots,\beta_\varphi\gamma_1,\gamma_2^\prime,\dots,\gamma_q^\prime}$$

Here $\Lambda, \Lambda_{t,t+1}$ are defined by (5), $\alpha_1\ne\beta_1,\dots,\alpha_\varphi\ne\beta_\varphi$, $\psi+\varphi+q=p$. In the case $\psi=q=0$ we have the previous relations (4).

We remark that if $\tau=0$ and $p=1$, or $p=2$, then the algebra $A_{m_1,\dots,m_p}(\Cal E,0)$ is the polynomial ring in variables $\{e_{\alpha_1,\dots,\alpha_p}\}$ over the field of meromorphic functions in variables $\{y_{\alpha,\beta}\}$. In the case $\tau=0$, $p>2$ the algebra $A_{m_1,\dots,m_p}(\Cal E,0)$ is commutative but it is not a polynomial ring, because the relations (6) take the form:
$$e_{\mu_1,\dots,\mu_{\psi-1},\mu_\psi,\alpha_1,\dots,\alpha_\varphi,\gamma_1,\gamma_2,\dots,\gamma_q}e_{\mu_1^\prime,\dots,\mu_{\psi-1}^\prime,\mu_\psi,\beta_1,\dots,\beta_\varphi,\gamma_1,\gamma_2^\prime,\dots,\gamma_q^\prime}=\eqno(7)$$
$$e_{\mu_1,\dots,\mu_{\psi-1},\mu_\psi,\beta_1,\dots,\beta_\varphi,\gamma_1,\gamma_2,\dots,\gamma_q}e_{\mu_1^\prime,\dots,\mu_{\psi-1}^\prime,\mu_\psi,\alpha_1,\dots,\alpha_\varphi,\gamma_1,\gamma_2^\prime,\dots,\gamma_q^\prime}$$

The algebra $A_{m_1,\dots,m_p}(\Cal E,\tau)$ is a flat deformation of the algebra of functions on the manifold defined by the equalities (7). It is easy to see that this manifold is rational and the general solution of the equalities (7) has a form: $e_{\alpha_1,\dots,\alpha_p}=e^{(1,2)}_{\alpha_1,\alpha_2}e^{(2,3)}_{\alpha_2,\alpha_3}\dots e^{(p-1,p)}_{\alpha_{p-1},\alpha_p}$, where $\{e^{(\varphi,\varphi+1)}_{\alpha_\varphi,\alpha_{\varphi+1}}\}$ are independent variables.

{\bf Proposition 2.} {\it There is the homomorphism of the algebras 

$$x: Q_{n,k}(\Cal E,\tau)\to A_{m_1,\dots,m_p}(\Cal E,\tau)$$
that acts on the generators of the algebra $Q_{n,k}(\Cal E,\tau)$ in the following way: $f(z_1,\dots,z_p)\in\Theta_{(n_1,\dots,n_p)}(\Gamma)$ is sent to} 
$$x(f)=\sum\Sb1\leqslant\alpha_1\leqslant m_1\\ \\\dots\dots\dots\dots\dots\dots\\ \\1\leqslant\alpha_p\leqslant m_p\endSb f(y_{1,\alpha_1},\dots,y_{p,\alpha_p})e_{\alpha_1,\dots,\alpha_p}$$

\newpage

\centerline{\bf \S3. Functional realization of the Poisson algebra $q_{n,k}\widehat\otimes q_{m,l}(\Cal E)$}
\medskip
Let $\frac{n}{k}=n_1-\frac{1}{n_2-\dots-\frac{1}{n_p}}$, $\frac{m}{l}=m_1-\frac{1}{m_2-\dots-\frac{1}{m_q}}$, where $n_1,\dots,n_p,m_1,\dots,m_q\geqslant2$. We denote $q_{n,k}\widehat\otimes q_{m,l}(\Cal E)=\botshave{\bigoplus_{\alpha,\beta\geqslant0}}P_{\alpha,\beta}$, where $P_{\alpha,\beta}$ is a space of meromorphic functions $f(x_{1,1},\dots,x_{p,1};\dots;x_{1,\alpha},\dots,x_{p,\alpha};y_{1,1},\dots,y_{q,1};\dots;y_{1,\beta},\dots,y_{q,\beta};z)$ in variables $\{x_{i,j}; 1\leqslant i\leqslant p, 1\leqslant j\leqslant\alpha\}$, $\{y_{i,j};1\leqslant i\leqslant q, 1\leqslant j\leqslant\beta\}$ and $z$ with the following properties:

1. Symmetry. For each $\sigma\in S_\alpha, \delta\in S_\beta$, $f$ is invariant with respect to permutations $x_{i,j}\mapsto x_{i,\sigma_j}, y_{i,j}\mapsto y_{i,\delta_j}$.

2. Periodicity and quasiperiodicity. $$f(x_{1,1},\dots,x_{t,\mu}+1,\dots,x_{p,\alpha};y_{1,1},\dots,y_{q,\beta};z)=f(x_{1,1},\dots,y_{q,\beta};z)$$

$$f(x_{1,1},\dots,x_{p,\alpha};y_{1,1},\dots,y_{t,\mu}+1,\dots,y_{q,\beta};z)=f(x_{1,1},\dots,y_{q,\beta};z)$$

$$f(x_{1,1},\dots,x_{t,\mu}+\eta,\dots,x_{p,\alpha};y_{1,1},\dots,y_{q,\beta};z)=e^{-2\pi i(n_tx_{t,\mu}-x_{t-1,\mu}-x_{t+1,\mu})}f(x_{1,1},\dots,y_{q,\beta},z)$$

$$f(x_{1,1},\dots,x_{p,\alpha};y_{1,1},\dots,y_{t,\mu}+\eta,\dots,y_{q,\beta};z)=e^{-2\pi i(m_ty_{t,\mu}-y_{t-1,\mu}-y_{t+1,\mu})}f(x_{1,1},\dots,y_{q,\beta};z)$$

Here $x_{0,\mu}=x_{p+1,\mu}=y_{0,\mu}=y_{q+1,\mu}=0$.

3. $f(x_{1,1},\dots,y_{q,\beta};z)$ as a function in variables $\{x_{t,\mu}, y_{t^\prime,\mu^\prime}\}$ is holomorphic outside the divisors $\{y_{1,\mu}-x_{p,\mu^\prime}-z=0; 1\leqslant\mu\leqslant\beta, 1\leqslant\mu^\prime\leqslant\alpha\}$ and has a pole of order $\leqslant1$ on these divisors. So $\hat f(x_{1,1},\dots,y_{q,\beta};z)=\bigg(\botshave{\prod\Sb1\leqslant\mu\leqslant\beta\\ \\1\leqslant\mu^\prime\leqslant\alpha\endSb}\theta(y_{1,\mu}-x_{p,\mu^\prime}-z)\bigg)f(x_{1,1},\dots,y_{q,\beta};z)$ as a function in variables $\{x_{t,\mu}, y_{t^\prime,\mu^\prime}\}$ is holomorphic.

4. Let $1\leqslant\mu_1\ne\mu_2,\mu_3\leqslant\beta, 1\leqslant\nu_1,\nu_2\ne\nu_3\leqslant\alpha$. Then $\hat f(x_{1,1},\dots,y_{q,\beta};z)=0$ on the affine subspaces of codimension 2 defined by the following relations: $y_{1,\mu_1}=y_{1,\mu_2}=x_{p,\nu_1}+z$ or $y_{1,\mu_3}=x_{p,\nu_2}+z=x_{p,\nu_3}+z$.

Particularly, the space $M=P_{0,0}$ is a field of meromorphic functions $f(z)$. It is possible to check, that the dimensions $\roman{dim}_MP_{\alpha,\beta}$ are finite for each $\alpha, \beta$ and the Hilbert function is 
$$\botshave{\sum_{\alpha, \beta\geqslant0}}\roman{dim}_MP_{\alpha,\beta}t_1^\alpha t_2^\beta=(1-t_1)^{-d(n_1,\dots,n_p)}(1-t_2)^{-d(m_1,\dots,m_q)}(1-t_1t_2)^{-d(n_1,\dots,n_p+m_1,\dots,m_q)}$$

We define the commutative associative product on the space $q_{n,k}\widehat\otimes q_{m,l}(\Cal E)$ in the following way: if $f\in P_{\alpha,\beta}, g\in P_{\alpha^\prime,\beta^\prime}$ then for $fg\in P_{\alpha+\alpha^\prime,\beta+\beta^\prime}$ we have:
$$fg(x_{1,1},\dots,x_{p,\alpha+\alpha^\prime};y_{1,1},\dots,y_{q,\beta+\beta^\prime};z)=\frac{1}{\alpha!\alpha^\prime!\beta!\beta^\prime!}\sum\Sb\sigma\in S_{\alpha+\alpha^\prime}\\ \\\delta\in S_{\beta+\beta^\prime}\endSb$$

$$f(x_{1,\sigma_1},\dots,x_{p,\sigma_\alpha};y_{1,\delta_1},\dots,y_{q,\delta_\beta};z)g(x_{1,\sigma_{\alpha+1}},\dots,x_{p,\sigma_{\alpha+\alpha^\prime}};y_{1,\delta_{\beta+1}},\dots,y_{q,\delta_{\beta+\beta^\prime}};z)$$

We define the Poisson bracket $\{,\}$ on the space $q_{n,k}\widehat\otimes q_{m,l}(\Cal E)$ in the following way: if $f\in P_{\alpha,0}, g\in P_{\beta,0}$ and $f, g$ do not depend on $z$, then the Poisson bracket $\{f,g\}$ is given by the formula (3). Similar formula gives the Poisson bracket $\{f,g\}$ in the case $f\in P_{0,\alpha}, g\in P_{0,\beta}$ and $f, g$ do not depend on $z$. Let $f\in P_{1,0}$, $g\in P_{0,1}$ and $f, g$ do not depend on $z$. In this case, for $\{f,g\}$ we have by definition:
$$\{f,g\}(x_1,\dots,x_p;y_1,\dots,y_q;z)=$$

$$f(x_1,\dots,x_p)\botshave{\sum_{1\leqslant t\leqslant q}}\frac{d(m_{t+1},\dots,m_q)}{d(m_1,\dots,m_q)}g^\prime_{y_t}(y_1,\dots,y_q)-$$

$$-g(y_1,\dots,y_q)\botshave{\sum_{1\leqslant t\leqslant p}}\frac{d(n_1,\dots,n_{t-1})}{d(n_1,\dots,n_p)}f^\prime_{x_t}(x_1,\dots,x_p)-$$

$$-\bigg(\frac{\theta^\prime(y_1-x_p-z)}{\theta(y_1-x_p-z)}-\pi i\bigg)f(x_1,\dots,x_p)g(y_1,\dots,y_q)$$

For $f\in P_{1,0}$ and $g\in P_{0,1}$ we assume:

$$\{f,z\}=\bigg(\frac{d(m_2,\dots,m_q)}{d(m_1,\dots,m_q)}+\frac{d(n_1,\dots,n_{p-1})+1}{d(n_1,\dots,n_p)}\bigg)f$$

$$\{g,z\}=-\bigg(\frac{1+d(m_2,\dots,m_q)}{d(m_1,\dots,m_q)}+\frac{d(n_1,\dots,n_{p-1})}{d(n_1,\dots,n_p)}\bigg)g$$

It is clear that using the Leibniz rule one can extend the Poisson bracket $\{,\}$ on the whole space $q_{n,k}\widehat\otimes q_{m,l}(\Cal E)$.

For description of this bracket let us consider the commutative associative algebra generated by all meromorphic functions in variables $\{x^{(\gamma)}_{\alpha,\beta}, y^{(\gamma^\prime)}_{\alpha^\prime,\beta^\prime},z; 1\leqslant\alpha\leqslant p, 1\leqslant\alpha^\prime\leqslant q; \beta, \beta^\prime, \gamma, \gamma^\prime\in \Bbb N\}$ and the elements $\{e(x^{(\gamma)}_{1,\beta},\dots,x^{(\gamma)}_{p,\beta}), e^\prime(y^{(\gamma^\prime)}_{1,\beta^\prime},\dots,y^{(\gamma^\prime)}_{q,\beta^\prime}); \beta, \beta^\prime, \gamma, \gamma^\prime\in\Bbb N\}$. We assume that if $u_\alpha=v_\alpha$ for some $\alpha$, then $e(u_1,\dots,u_p)e(v_1,\dots,v_p)=0, e^\prime(u_1,\dots,u_q)e^\prime(v_1,\dots,v_q)=0$.

We define the Poisson structure on this algebra in the following way:

$$\{x^{(\gamma)}_{\alpha,\beta},x^{(\gamma^\prime)}_{\alpha^\prime,\beta^\prime}\}=
\{x^{(\gamma)}_{\alpha,\beta},y^{(\gamma^\prime)}_{\alpha^\prime,\beta^\prime}\}=
\{y^{(\gamma)}_{\alpha,\beta},y^{(\gamma^\prime)}_{\alpha^\prime,\beta^\prime}\}=\{x^{(\gamma)}_{\alpha,\beta},z\}=\{y^{(\gamma^\prime)}_{\alpha^\prime,\beta^\prime},z\}=0$$

$$\{e(u_1,\dots,u_p),e(v_1,\dots,v_p)\}=\bigg(\frac{\theta^\prime(v_1-u_1)}{\theta(v_1-u_1)}+\frac{\theta^\prime(v_p-u_p)}{\theta(v_p-u_p)}-2\pi i\bigg)e(u_1,\dots,u_p)e(v_1,\dots,v_p)+$$

$$\theta^\prime(0)\botshave{\sum_{1\leqslant\alpha< p}}\frac{\theta(v_\alpha+v_{\alpha+1}-u_\alpha-u_{\alpha+1})}{\theta(v_\alpha-u_\alpha)\theta(v_{\alpha+1}-u_{\alpha+1})}e(v_1,\dots,v_\alpha,u_{\alpha+1},\dots,u_p)e(u_1,\dots,u_\alpha,v_{\alpha+1},\dots,v_p)$$

$$\{e^\prime(u_1,\dots,u_q),e^\prime(v_1,\dots,v_q)\}=\bigg(\frac{\theta^\prime(v_1-u_1)}{\theta(v_1-u_1)}+\frac{\theta^\prime(v_q-u_q)}{\theta(v_q-u_q)}-2\pi i\bigg)e^\prime(u_1,\dots,u_q)e^\prime(v_1,\dots,v_q)+$$

$$\theta^\prime(0)\botshave{\sum_{1\leqslant\alpha<q}}\frac{\theta(v_\alpha+v_{\alpha+1}-u_\alpha-u_{\alpha+1})}{\theta(v_\alpha-u_\alpha)\theta(v_{\alpha+1}-u_{\alpha+1})}e^\prime(v_1,\dots,v_\alpha,u_{\alpha+1},\dots,u_q)e^\prime(u_1,\dots,u_\alpha,v_{\alpha+1},\dots,v_q)$$

$$\{e(u_1,\dots,u_p),e^\prime(v_1,\dots,v_q)\}=\bigg(\frac{\theta^\prime(u_p-v_1+z)}{\theta(u_p-v_1+z)}-\pi i\bigg)e(u_1,\dots,u_p)e^\prime(v_1,\dots,v_q)$$

$$\{e(u_1,\dots,u_p),x^{(\gamma)}_{\alpha,\beta}\}=-\frac{d(n_1,\dots,n_{\alpha-1})+d(n_{\alpha+1},\dots,n_p)}{d(n_1,\dots,n_p)}e(u_1,\dots,u_p)$$

$$\{e(u_1,\dots,u_p),y^{(\gamma)}_{\alpha,\beta}\}=\frac{d(m_{\alpha+1},\dots,m_q)}{d(m_1,\dots,m_q)}e(u_1,\dots,u_p)$$

$$\{e(u_1,\dots,u_p),z\}=\bigg(\frac{d(m_2,\dots,m_q)}{d(m_1,\dots,m_q)}+\frac{d(n_1,\dots,n_{p-1})+1}{d(n_1,\dots,n_p)}\bigg)e(u_1,\dots,u_p)$$

$$\{e^\prime(u_1,\dots,u_q),y^{(\gamma)}_{\alpha,\beta}\}=-\frac{d(m_1,\dots,m_{\alpha-1})+d(m_{\alpha+1},\dots,m_q)}{d(m_1,\dots,m_q)}e^\prime(u_1,\dots,u_q)$$

$$\{e^\prime(u_1,\dots,u_q),x^{(\gamma)}_{\alpha,\beta}\}=\frac{d(n_1,\dots,n_{\alpha-1})}{d(n_1,\dots,n_p)}e^\prime(u_1,\dots,u_q)$$

$$\{e^\prime(u_1,\dots,u_q),z\}=-\bigg(\frac{1+d(m_2,\dots,m_q)}{d(m_1,\dots,m_q)}+\frac{d(n_1,\dots,n_{p-1})}{d(n_1,\dots,n_p)}\bigg)e^\prime(u_1,\dots,u_q)$$

For $f\in P_{\alpha,\beta}$ we define the element $X_f$ in the following way:

$$X_f=\botshave{\sum_{\gamma_{\mu,\nu},\gamma^\prime_{\mu^\prime,\nu^\prime}\in\Bbb N}}f(x^{(\gamma_{1,1})}_{1,1},\dots,x^{(\gamma_{p,\alpha})}_{p,\alpha};y^{(\gamma^\prime_{1,1})}_{1,1},\dots,y^{(\gamma^\prime_{q,\beta})}_{q,\beta};z)\times$$

$$\times e(x_{1,1}^{(\gamma_{1,1})},\dots,x_{p,1}^{(\gamma_{p,1})})\dots e(x_{1,\alpha}^{(\gamma_{1,\alpha})},\dots,x_{p,\alpha}^{(\gamma_{p,\alpha})})e^\prime(y_{1,1}^{(\gamma^\prime_{1,1})},\dots,y_{q,1}^{(\gamma^\prime_{q,1})})\dots e^\prime(y_{1,\beta}^{(\gamma^\prime_{1,\beta})},\dots,y_{q,\beta}^{(\gamma^\prime_{q,\beta})})$$

The Poisson bracket and the product on the space $q_{n,k}\widehat\otimes q_{m,l}(\Cal E)$ are defined by the following formulas: $$X_{\{f,g\}}=\{X_f,X_g\}, X_{fg}=X_fX_g$$

\newpage

\centerline{\bf \S4. Functional realization of the Poisson algebra $q_{n_1,k_1}\widehat\otimes\dots\widehat\otimes q_{n_h,k_h}(\Cal E)$}
\medskip

Let $\frac{n_\alpha}{k_\alpha}=n_{1,\alpha}-\frac{1}{n_{2,\alpha}-\dots-\frac{1}{n_{p_\alpha,\alpha}}}$, where $1\leqslant\alpha\leqslant h$; $n_{\nu,\alpha}\geqslant2$. Let $q_{n_1,k_1}\widehat\otimes\dots\widehat\otimes q_{n_h,k_h}(\Cal E)=\botshave{\bigoplus_{\alpha_1,\dots,\alpha_h\geqslant0}}P_{\alpha_1,\dots,\alpha_h}$, where $P_{\alpha_1,\dots,\alpha_h}$ is the space of meromorphic functions in variables $\{x_{\mu,\lambda,t}, z_{\nu,\nu+1}; 1\leqslant t\leqslant h, 1\leqslant\mu\leqslant p_t, 1\leqslant\lambda\leqslant\alpha_t, 1\leqslant\nu\leqslant h-1\}$ with the following properties:

1.Symmetry. Let $\sigma_t\in S_{\alpha_t}$ for $1\leqslant t\leqslant h$. Then $f$ is invariant with respect to permutations: $x_{\mu,\lambda,t}\mapsto x_{\mu,\sigma_t(\lambda),t}$.

2. Periodicity and quasiperiodicity. $$f(x_{1,1,1},\dots,x_{\mu,\lambda,t}+1,\dots,x_{p_h,\alpha_h,h};z_{1,2},\dots,z_{h-1,h})=f(x_{1,1,1},\dots,z_{h-1,h})$$

$$f(x_{1,1,1},\dots,x_{\mu,\lambda,t}+\eta,\dots,x_{p_h,\alpha_h,h};z_{1,2},\dots,z_{h-1,h})=$$

$$e^{-2\pi i(n_{\mu,t}x_{\mu,\lambda,t}-x_{\mu-1,\lambda,t}-x_{\mu+1,\lambda,t})}f(x_{1,1,1},\dots,z_{h-1,h})$$

3. $f$ as a function in variables $\{x_{\mu,\lambda,t}\}$ is holomorphic outside the divisors $x_{1,\mu,t+1}-x_{p_t,\mu^\prime,t}-z_{t,t+1}=0$ and has a pole of order $\leqslant1$ on these divisors. So the function $\hat f=\bigg(\botshave{\prod\Sb1\leqslant\mu\leqslant\alpha_{t+1}\\ \\1\leqslant\mu^\prime\leqslant\alpha_t\\ \\1\leqslant t<h\endSb}\theta(x_{1,\mu,t+1}-x_{p_t,\mu^\prime,t}-z_{t,t+1})\bigg)f$ as a function in variables $\{x_{\mu,\lambda,t}\}$ is holomorphic.

4. $\hat f=0$ on the affine subspaces of codimension 2 defined by the relations $x_{1,\mu_1,t+1}=x_{1,\mu_2,t+1}=x_{p_t,\mu_3,t}+z_{t,t+1}$ or $x_{1,\mu_1,t+1}=x_{p_t,\mu_2,t}+z_{t,t+1}=x_{p_t,\mu_3,t}+z_{t,t+1}$.

Particularly, the space $M=P_{0.\dots,0}$ is the field of meromorphic functions in variables $z_{1,2},\dots,z_{h-1,h}$. It is possible to check that the dimensions $\roman{dim}_MP_{\alpha_1,\dots,\alpha_h}$ are finite and the Hilbert function is 
$$\botshave{\sum_{\alpha_1,\dots,\alpha_h\geqslant0}}\roman{dim}_MP_{\alpha_1,\dots,\alpha_h}t_1^{\alpha_1}\dots t_h^{\alpha_h}=\botshave{\prod_{1\leqslant\lambda\leqslant\nu\leqslant h}}(1-t_\lambda t_{\lambda+1}\dots t_\nu)^{-d(N_\lambda\Delta N_{\lambda+1}\Delta\dots\Delta N_\nu)}\eqno(8)$$

Here $N_t=(n_{1,t},\dots,n_{p_t,t})$ for $1\leqslant t\leqslant h$. See {\bf Notations} for definition $\Delta$.

For definition of the product and the Poisson bracket on the space $q_{n_1,k_1}\widehat\otimes\dots\widehat\otimes q_{n_h,k_h}(\Cal E)$ we consider the commutative associative algebra generated by all meromorphic functions in variables $\{x^{(\gamma)}_{\mu,\lambda,t}, z_{\nu,\nu+1}; 1\leqslant\mu\leqslant p_t, 1\leqslant t\leqslant h, 1\leqslant\nu\leqslant h-1, \lambda,\gamma\in\Bbb N\}$ and the elements $e_t(x^{(\gamma)}_{1,\lambda,t},\dots,x^{(\gamma)}_{p_t,\lambda,t})$. We assume that if $u_\alpha=v_\alpha$ for some $1\leqslant\alpha\leqslant p_t$, then $e_t(u_1,\dots,u_{p_t})e_t(v_1,\dots,v_{p_t})=0$. We define the Poisson bracket on this algebra in the following way: 

$$\{x^{(\gamma)}_{\mu,\lambda,t},x^{(\gamma^\prime)}_{\mu^\prime,\lambda^\prime,t^\prime}\}=\{x^{(\gamma)}_{\mu,\lambda,t},z_{\nu,\nu+1}\}=\{z_{\nu,\nu+1},z_{\nu^\prime,\nu^\prime+1}\}=0$$

$$\{e_t(u_1,\dots,u_{p_t}),x^{(\gamma)}_{\mu,\lambda,t}\}=-\frac{d(n_{1,t},\dots,n_{\mu-1,t})+d(n_{\mu+1,t},\dots,n_{p_t,t})}{d(n_{1,t},\dots,n_{p_t,t})}e_t(u_1,\dots,u_{p_t})$$

$$\{e_t(u_1,\dots,u_{p_t}),x^{(\gamma)}_{\mu,\lambda,t+1}\}=\frac{d(n_{\mu+1,t+1},\dots,n_{p_{t+1},t+1})}{d(n_{1,t+1},\dots,n_{p_{t+1},t+1})}e_t(u_1,\dots,u_{p_t})$$

$$\{e_{t+1}(u_1,\dots,u_{p_{t+1}}),x^{(\gamma)}_{\mu,\lambda,t}\}=\frac{d(n_{1,t},\dots,n_{\mu-1,t})}{d(n_{1,t},\dots,n_{p_t,t})}e_{t+1}(u_1,\dots,u_{p_{t+1}})$$

$$\{e_t(u_1,\dots,u_{p_t}),x^{(\gamma)}_{\mu,\lambda,t^\prime}\}=0 \text{, here } |t-t^\prime|>1$$

$$\{e_t(u_1,\dots,u_{p_t}),z_{t,t+1}\}=$$

$$\bigg(\frac{d(n_{2,t+1},\dots,n_{p_{t+1},t+1})}{d(n_{1,t+1},\dots,n_{p_{t+1},t+1})}+\frac{d(n_{1,t},\dots,n_{p_t-1,t})+1}{d(n_{1,t},\dots,n_{p_t,t})}\bigg)e_t(u_1,\dots,u_{p_t})$$

$$\{e_{t+1}(u_1,\dots,u_{p_{t+1}}),z_{t,t+1}\}=$$

$$-\bigg(\frac{d(n_{2,t+1},\dots,n_{p_{t+1},t+1})+1}{d(n_{1,t+1},\dots,n_{p_{t+1},t+1})}+\frac{d(n_{1,t},\dots,n_{p_t-1,t})}{d(n_{1,t},\dots,n_{p_t,t})}\bigg)e_{t+1}(u_1,\dots,u_{p_{t+1}})$$

$$\{e_t(u_1,\dots,u_{p_t}),z_{t^\prime,t^\prime+1}\}=0 \text{, here }t\ne t^\prime, t\ne t^\prime+1$$ 

$$\{e_t(u_1,\dots,u_{p_t}),e_t(v_1,\dots,v_{p_t})\}=$$

$$\bigg(\frac{\theta^\prime(v_1-u_1)}{\theta(v_1-u_1)}+\frac{\theta^\prime(v_{p_t}-u_{p_t})}{\theta(v_{p_t}-u_{p_t})}-2\pi i\bigg)e_t(u_1,\dots,u_{p_t})e_t(v_1,\dots,v_{p_t})+$$

$$\theta^\prime(0)\botshave{\sum_{1\leqslant\alpha<p_t}}\frac{\theta(v_\alpha+v_{\alpha+1}-u_\alpha-u_{\alpha+1})}{\theta(v_\alpha-u_\alpha)\theta(v_{\alpha+1}-u_{\alpha+1})}e_t(v_1,\dots,v_\alpha,u_{\alpha+1},\dots,u_{p_t})e(u_1,\dots,u_\alpha,v_{\alpha+1},\dots,v_{p_t})$$

$$\{e_t(u_1,\dots,u_{p_t}),e_{t+1}(v_1,\dots,v_{p_{t+1}})\}=$$

$$\bigg(\frac{\theta^\prime(u_{p_t}-v_1+z_{t,t+1})}{\theta(u_{p_t}-v_1+z_{t,t+1})}-\pi i\bigg)e_t(u_1,\dots,u_{p_t})e_{t+1}(v_1,\dots,v_{p_{t+1}})$$

$$\{e_t(u_1,\dots,u_{p_t}),e_{t^\prime}(v_1,\dots,v_{p_{t^\prime}})=0\text{, here }|t-t^\prime|>1$$

For $f\in P_{\alpha_1,\dots,\alpha_h}$ we define the element $X_f$ by the following formula:

$$X_f=\botshave{\sum_{\gamma_{\mu,\lambda,\nu\in \Bbb N}}}f(x^{(\gamma_{1,1,1})}_{1,1,1},\dots,x^{(\gamma_{p_h,\alpha_h,h})}_{p_h,\alpha_h,h},z_{1,2},\dots,z_{h-1,h})\prod\Sb1\leqslant t\leqslant h\\ \\1\leqslant\lambda\leqslant\alpha_t\endSb e_t(x^{(\gamma_{1,\lambda,t})}_{1,\lambda,t},\dots,x^{(\gamma_{p_t,\lambda,t})}_{p_t,\lambda,t})$$

{\bf Proposition 3.} {\it The following formulas define the product and the Poisson bracket on the space $q_{n_1,k_1}\widehat\otimes\dots\widehat\otimes q_{n_h,k_h}(\Cal E)$:}

$$X_{fg}=X_fX_g, X_{\{f,g\}}=\{X_f,X_g\}$$

\newpage
\centerline{\bf \S5. Bosonization of the algebra $Q_{n_1,k_1}\widehat\otimes\dots\widehat\otimes Q_{n_h,k_h}(\Cal E,\tau)$}

\medskip

Let $A_{l_{1,1},\dots,l_{p_1,1};\dots;l_{1,h},\dots,l_{p_h,h}}(\Cal E,\tau)$ be the associative algebra generated by all meromorphic functions in variables $\{y_{\lambda,\mu,\gamma}, z_{\gamma^\prime,\gamma^\prime+1}; 1\leqslant\gamma\leqslant h, 1\leqslant\gamma^\prime\leqslant h-1, 1\leqslant\lambda\leqslant p_\gamma, 1\leqslant\mu\leqslant l_{\lambda,\gamma}\}$ and the elements $\{e^{(\gamma)}_{\alpha_1,\dots,\alpha_{p_\gamma}}; 1\leqslant\gamma\leqslant h, 1\leqslant\alpha_\lambda\leqslant l_{\lambda,\gamma},1\leqslant\lambda\leqslant p_\gamma\}$. We assume that the following relations hold:

$$[y_{\lambda,\mu,\gamma},y_{\lambda^\prime,\mu^\prime,\gamma^\prime}]=[y_{\lambda,\mu,\gamma},z_{t,t+1}]=[z_{t,t+1},z_{t^\prime,t^\prime+1}]=0$$

$$e^{(\gamma)}_{\alpha_1,\dots,\alpha_{p_\gamma}}y_{\lambda,\beta,\gamma}=\bigg(y_{\lambda,\beta,\gamma}-\frac{d(n_{1,\gamma},\dots,n_{\lambda-1,\gamma})+d(n_{\lambda+1,\gamma},\dots,n_{p_\gamma,\gamma})}{d(n_{1,\gamma},\dots,n_{p_\gamma,\gamma})}\tau\bigg)e^{(\gamma)}_{\alpha_1,\dots,\alpha_{p_\gamma}}$$
Here $\beta\ne\alpha_\lambda$.

$$e^{(\gamma)}_{\alpha_1,\dots,\alpha_{p_\gamma}}y_{\lambda,\alpha_\lambda,\gamma}=\bigg(y_{\lambda,\alpha_\lambda,\gamma}+\bigg(1-\frac{d(n_{1,\gamma},\dots,n_{\lambda-1,\gamma})+d(n_{\lambda+1,\gamma},\dots,n_{p_\gamma,\gamma})}{d(n_{1,\gamma},\dots,n_{p_\gamma,\gamma})}\bigg)\tau\bigg)e^{(\gamma)}_{\alpha_1,\dots,\alpha_{p_\gamma}}$$

$$e^{(\gamma)}_{\alpha_1,\dots,\alpha_{p_\gamma}}y_{\lambda,\beta,\gamma+1}=\bigg(y_{\lambda,\beta,\gamma+1}+\frac{d(n_{\lambda+1,\gamma+1},\dots,n_{p_{\gamma+1},\gamma+1})}{d(n_{1,\gamma+1},\dots,n_{p_{\gamma+1},\gamma+1})}\tau\bigg)e^{(\gamma)}_{\alpha_1,\dots,\alpha_{p_\gamma}}$$

$$e^{(\gamma+1)}_{\alpha_1,\dots,\alpha_{p_{\gamma+1}}}y_{\lambda,\beta,\gamma}=\bigg(y_{\lambda,\beta,\gamma}+\frac{d(n_{1,\gamma},\dots,n_{\lambda-1,\gamma})}{d(n_{1,\gamma},\dots,n_{p_\gamma,\gamma})}\tau\bigg)e^{(\gamma+1)}_{\alpha_1,\dots,\alpha_{p_{\gamma+1}}}$$

$$e^{(\gamma)}_{\alpha_1,\dots,\alpha_{p_\gamma}}y_{\lambda,\beta,\gamma^\prime}=y_{\lambda,\beta,\gamma^\prime}e^{(\gamma)}_{\alpha_1,\dots,\alpha_{p_\gamma}}$$
Here $|\gamma-\gamma^\prime|>1$.

$$e^{(\gamma)}_{\alpha_1,\dots,\alpha_{p_\gamma}}z_{\gamma,\gamma+1}=$$

$$\bigg(z_{\gamma,\gamma+1}+\bigg(\frac{d(n_{2,\gamma+1},\dots,n_{p_{\gamma+1},\gamma+1})}{d(n_{1,\gamma+1},\dots,n_{p_{\gamma+1},\gamma+1})}+\frac{d(n_{1,\gamma},\dots,n_{p_\gamma-1,\gamma})+1}{d(n_{1,\gamma},\dots,n_{p_\gamma,\gamma})}\bigg)\tau\bigg)e^{(\gamma)}_{\alpha_1,\dots,\alpha_{p_\gamma}}$$

$$e^{(\gamma+1)}_{\alpha_1,\dots,\alpha_{p_{\gamma+1}}}z_{\gamma,\gamma+1}=$$

$$\bigg(z_{\gamma,\gamma+1}-\bigg(\frac{d(n_{2,\gamma+1},\dots,n_{p_{\gamma+1},\gamma+1})+1}{d(n_{1,\gamma+1},\dots,n_{p_{\gamma+1},\gamma+1})}+\frac{d(n_{1,\gamma},\dots,n_{p_\gamma-1,\gamma})}{d(n_{1,\gamma},\dots,n_{p_\gamma,\gamma})}\bigg)\tau\bigg)e^{(\gamma+1)}_{\alpha_1,\dots,\alpha_{p_{\gamma+1}}}$$

$$e^{(\gamma)}_{\alpha_1,\dots,\alpha_{p_\gamma}}z_{\gamma^\prime,\gamma^\prime+1}=z_{\gamma^\prime,\gamma^\prime+1}e^{(\gamma)}_{\alpha_1,\dots,\alpha_{p_\gamma}}$$
Here $\gamma\ne\gamma^\prime$ and $\gamma\ne\gamma^\prime+1$.

$$e^{(\gamma)}_{\alpha_1,\dots,\alpha_{p_\gamma}}e^{(\gamma^\prime)}_{\beta_1,\dots,\beta_{p_{\gamma^\prime}}}=e^{(\gamma^\prime)}_{\beta_1,\dots,\beta_{p_{\gamma^\prime}}}e^{(\gamma)}_{\alpha_1,\dots,\alpha_{p_\gamma}}$$
Here $|\gamma-\gamma^\prime|>1$.

$$e^{(\gamma)}_{\alpha_1,\dots,\alpha_{p_\gamma}}e^{(\gamma+1)}_{\beta_1,\dots,\beta_{p_{\gamma+1}}}=$$

$$-e^{2\pi i(y_{1,\beta_1,\gamma+1}-y_{p_\gamma,\alpha_{p_\gamma},\gamma}-z_{\gamma,\gamma+1})}\frac{\theta(y_{p_\gamma,\alpha_{p_\gamma},\gamma}-y_{1,\beta_1,\gamma+1}+z_{\gamma,\gamma+1}+\frac{1}{2}\tau)}{\theta(y_{1,\beta_1,\gamma+1}-y_{p_\gamma,\alpha_{p_\gamma},\gamma}-z_{\gamma,\gamma+1}+\frac{1}{2}\tau)}e^{(\gamma+1)}_{\beta_1,\dots,\beta_{p_{\gamma+1}}}e^{(\gamma)}_{\alpha_1,\dots,\alpha_{p_\gamma}}$$

$$e^{(\gamma)}_{\alpha_1,\dots,\alpha_{p_\gamma}}e^{(\gamma)}_{\beta_1,\dots,\beta_{p_\gamma}}=\Phi e^{(\gamma)}_{\beta_1,\dots,\beta_{p_\gamma}}e^{(\gamma)}_{\alpha_1,\dots,\alpha_{p_\gamma}}+\eqno(9)$$

$$+\botshave{\sum_{1\leqslant t\leqslant p_\gamma-1}}\Phi_{t,t+1}e^{(\gamma)}_{\beta_1,\dots,\beta_t,\alpha_{t+1},\dots,\alpha_{p_\gamma}}e^{(\gamma)}_{\alpha_1,\dots,\alpha_t,\beta_{t+1},\dots,\beta_{p_\gamma}}$$

Here $\alpha_1\ne\beta_1,\dots,\alpha_{p_\gamma}\ne\beta_{p_\gamma}$ and

$$\Phi=\frac{e^{-2\pi i\tau}\theta(y_{1,\beta_1,\gamma}-y_{1,\alpha_1,\gamma})\theta(y_{p_\gamma,\beta_{p_\gamma},\gamma}-y_{p_\gamma,\alpha_{p_\gamma},\gamma}+\tau)}{\theta(y_{1,\beta_1,\gamma}-y_{1,\alpha_1,\gamma}-\tau)\theta(y_{p_\gamma,\beta_{p_\gamma},\gamma}-y_{p_\gamma,\alpha_{p_\gamma},\gamma})}\eqno(10)$$

$$\Phi_{t,t+1}=\frac{e^{-2\pi i\tau}\theta(\tau)\theta(y_{1,\beta_1,\gamma}-y_{1,\alpha_1,\gamma})}{\theta(y_{1,\beta_1,\gamma}-y_{1,\alpha_1,\gamma}-\tau)}\cdot\frac{\theta(y_{t,\beta_t,\gamma}+y_{t+1,\beta_{t+1},\gamma}-y_{t,\alpha_t,\gamma}-y_{t+1,\alpha_{t+1},\gamma})}{\theta(y_{t,\beta_t,\gamma}-y_{t,\alpha_t,\gamma})\theta(y_{t+1,\beta_{t+1},\gamma}-y_{t+1,\alpha_{t+1},\gamma})}$$

In the general case, if some indexes are the same, then we have:

$$e^{(\gamma)}_{\mu_1,\dots,\mu_{\psi-1},\mu_\psi,\alpha_1,\dots,\alpha_\varphi,\gamma_1,\gamma_2,\dots,\gamma_q}e^{(\gamma)}_{\mu_1^\prime,\dots,\mu_{\psi-1}^\prime,\mu_\psi,\beta_1,\dots,\beta_\varphi,\gamma_1,\gamma_2^\prime,\dots,\gamma_q^\prime}=$$

$$\Phi e^{(\gamma)}_{\mu_1,\dots,\mu_{\psi-1},\mu_\psi,\beta_1,\dots,\beta_\varphi,\gamma_1,\gamma_2,\dots,\gamma_q}e^{(\gamma)}_{\mu_1^\prime,\dots,\mu_{\psi-1}^\prime,\mu_\psi,\alpha_1,\dots,\alpha_\varphi,\gamma_1,\gamma_2^\prime,\dots,\gamma_q^\prime}+$$

$$\botshave{\sum_{1\leqslant t<\varphi}}\Phi_{t,t+1}e^{(\gamma)}_{\mu_1,\dots,\mu_{\psi-1},\mu_\psi,\beta_1,\dots,\beta_t,\alpha_{t+1},\dots,\alpha_\varphi\gamma_1,\gamma_2,\dots,\gamma_q}e^{(\gamma)}_{\mu_1^\prime,\dots,\mu_{\psi-1}^\prime,\mu_\psi,\alpha_1,\dots,\alpha_t,\beta_{t+1},\dots,\beta_\varphi\gamma_1,\gamma_2^\prime,\dots,\gamma_q^\prime}$$
Here $\Phi, \Phi_{t,t+1}$ are defined by (10), $\alpha_1\ne\beta_1,\dots,\alpha_\varphi\ne\beta_\varphi$, $\psi+\varphi+q=p_\gamma$. In the case $\psi=q=0$ we have the previous relations (9).

{\bf Proposition 4.} {\it There is a family of the associative algebras  $Q_{n_1,k_1}\widehat\otimes\dots\widehat\otimes Q_{n_h,k_h}(\Cal E,\tau)$ that is the quantization of the Poisson algebra $q_{n_1,k_1}\widehat\otimes\dots\widehat\otimes q_{n_h,k_h}(\Cal E)$ and the homomorphism of the associative algebras $$x:Q_{n_1,k_1}\widehat\otimes\dots\widehat\otimes Q_{n_h,k_h}(\Cal E,\tau)\to P_{l_{1,1},\dots,l_{p_h,h}}(\Cal E,\tau)$$ such that for the element $f\in P_{0,\dots,0,1,0,\dots,0}$ (here 1 is on $\gamma$-th place) we have:

$$x(f)=\sum\Sb1\leqslant\alpha_1\leqslant l_{1,\gamma}\\ \\\dots\dots\dots\dots\dots\dots\\ \\1\leqslant\alpha_{p_\gamma}\leqslant l_{p_\gamma,\gamma}\endSb f(y_{1,\alpha_1,\gamma},\dots,y_{p_\gamma,\alpha_{p_\gamma},\gamma};z_{1,2},\dots,z_{h-1,h})e^{(\gamma)}_{\alpha_1,\dots,\alpha_{p_\gamma}}$$

The algebras $Q_{n_1,k_1}\widehat\otimes\dots\widehat\otimes Q_{n_h,k_h}(\Cal E,\tau)$ and $q_{n_1,k_1}\widehat\otimes\dots\widehat\otimes q_{n_h,k_h}(\Cal E)$ have the same Hilbert function (see (8)).}

\head Acknowledgments\endhead

This work was supported by the grant INTAS--OPEN--97--1312.

\newpage
\centerline{\bf References}
\medskip

1. A.V.Odesskii and B.L.Feigin, "Sklyanin's elliptic algebras", Funkts. Anal. Philozhen., 23, No. 3, 45-54 (1989).

2. A.V.Odesskii and B.L.Feigin, "Constructions of Sklyanin elliptic algebras and quantum $R$-matrices", Funkts. Anal. Philozen. 27, No. 1, 37-45 (1993).

3. A.V.Odesskii and B.L.Feigin, "Elliptic deformation of current algebras and their representations by difference operators." Funct. Anal. Appl. 31 (1997), No 3, 193-203 (1998).

4. A.V.Odesskii and B.L.Feigin, "Quantized moduli spaces of the bundles on the elliptic curve and their applications.", MPI Preprint Series 1998 (132), Bonn, Germany. (math.QA/9812059)

5. B.L.Feigin and A.V.Odesskii, Vector Bundles on Elliptic Curve and Sklyanin Algebras, RIMS-1032, September 1995, Kyoto University, Kyoto, Japan. (q-alg/9509021)

\enddocument